\documentclass[12pt]{amsart}

\advance\hoffset-15truemm

\let\Cal\mathcal
\let\Bbb\mathbb
\let\frak\mathfrak
\let\phi\varphi

\newcommand{\x}{\times}
\renewcommand{\o}{\circ}
\newcommand{\al}{\alpha}
\newcommand{\be}{\beta}
\newcommand{\Ga}{\Gamma}

\newcommand{\om}{\omega}
\newcommand{\Om}{\Omega}
\newcommand{\la}{\lambda}
\newcommand{\ph}{\phi}
\newcommand{\Ph}{\Phi}
\newcommand{\ps}{\psi}
\newcommand{\ka}{\kappa}

\newcommand{\ze}{\zeta}

\newcommand{\tcg}{{\tilde{\Cal G}}}

\newcommand{\Ad}{\operatorname{Ad}}

\newcommand{\fg}{{\mathfrak g}}

\newtheorem*{prop*}{Proposition}

\newtheorem*{thm*}{Theorem}

\newtheorem*{lem*}{Lemma}

\newtheorem*{kor*}{Corollary}

\theoremstyle{definition}

\newtheorem*{definition*}{Definition}

\newtheorem*{example*}{Example}

\begin{document}

\title{On left invariant CR structures on $SU(2)$}
\author{Andreas \v Cap}
\thanks{The author was supported by project P15747-N05 of the ``Fonds zur
F\"orderung der wissenschaftlichen Forschung'' (FWF) }
\date{March 2006}

\address{Institut f\"ur Mathematik, Universit\"at Wien,
Nordbergstra\ss e~15, A--1090 Wien, Austria, and International Erwin
Schr\"odinger Institute for Mathematical Physics, Boltzmanngasse 9,
A-1090 Wien, Austria}

\email{andreas.cap@esi.ac.at}

\begin{abstract}
There is a well known one--parameter family of left invariant CR
structures on $SU(2)\cong S^3$. We show how purely algebraic methods
can be used to explicitly compute the canonical Cartan connections
associated to these structures and their curvatures. We also obtain
explicit descriptions of tractor bundles and tractor connections. 
\end{abstract}


\maketitle

\section{Introduction}\label{1}
Three dimensional CR structures are among the examples of geometric
structures for which Elie Cartan constructed an associated normal
Cartan connection, see \cite{Cartan}. The homogeneous model for this
geometry is $S^3$, viewed as a quotient of the semisimple group
$G:=PSU(2,1)$ by a parabolic subgroup $P$, so three dimensional CR
structures form an example of a \textit{parabolic geometry}.

This example is remarkable in many respects. On the one hand, it is
sufficiently complicated to incorporate many of the features of
general parabolic geometries. On the other hand, the low dimension of
the group $G$ and the fact that all important natural bundles are
(either real or complex) line bundles, and hence all sections can
locally be viewed as functions, simplify matters considerably. In
fact, Cartan was even able to describe an algorithm for computing the
essential curvature invariant of such structures. Moreover, many
questions that have to be attacked using representation theory for
general parabolic geometries can be easily solved directly in this
case. An example for this is provided by the analysis of possible
dimensions of automorphism groups in \cite{Srni04}.

Returning to the homogeneous model $S^3=G/P$, consider the compact
subgroup $K=SU(2)\subset G$. Acting with elements of $K$ induces a
diffeomorphism $K\cong S^3$, so we can actually view the standard CR
structure on $S^3$ as a left invariant structure on $K$.  In this
picture, the structure can be easily obtained from data on the Lie
algebra $\frak k$ of $K$. These data admit an evident one--parameter
deformation, which gives rise to a one parameter family of left
invariant CR structures on $K$. The aim of this article is to show
that the canonical Cartan connections associated to these CR
structures can be computed using only linear algebra. On the way, one
also gets explicit formulae for their curvatures. Finally, we also
describe all tractor bundles and normal tractor connections
explicitly. These developments should also serve as a basis for a more
detailed analysis of these structures and as a prototype for dealing
with general left invariant parabolic geometries.

\noindent
\textbf{Acknowledgment}: I would like to thank Olivier Biquard for
bringing this example to my attention during a discussion at the
Winter School in Srn\'\i{}.

\section{Left invariant CR structures on $SU(2)$}\label{2}
\subsection{3--dimensional CR structures}\label{2.1}
Recall that a CR structure on a $3$--manifold $M$ is given by a
complex line subbundle $H\subset TM$, which defines a contact
structure on $M$. The subbundle $H\subset TM$ is called the \textit{CR
  subbundle}. Equivalently, we have a rank two subbundle $H\subset TM$
endowed with a complex structure $J:H\to H$ such that for one (or
equivalently any) locally non vanishing section $\xi\in\Ga(H)$ the
vector fields $\xi$, $J(\xi)$ and $[\xi,J(\xi)]$ form a local frame
for $TM$. In contrast to higher dimensions, there is no condition of
partial integrability or integrability in dimension $3$.

Given two CR structures, there is an evident notion of a (local) CR
diffeomorphism. This is a (local) diffeomorphism $f$, such that for
each point $x$ the tangent map $T_xf$ maps the CR subbundle to the CR
subbundle and the restriction of $T_xf$ to the CR subbundle is complex
linear. 

The basic examples for such structures are provided by generic real
hypersurfaces in two dimensional complex manifolds. If $(\tilde
M,\tilde J)$ is a two dimensional complex manifold and $M\subset\tilde
M$ is a real hypersurface, then for each $x\in M$ the tangent space
$T_xM$ has real dimension $3$ and sits in $T_x\tilde M$, which is a
two dimensional complex vector space. Now $H_x:=T_xM\cap\tilde
J(T_xM)$ is a complex subspace of $T_x\tilde M$, which evidently must
have complex dimension one. By construction, the spaces $H_x$ fit
together to define a complex line subbundle $H\subset TM$, with the
complex structure $J$ given by the restriction of $\tilde J$.
Generically, the subbundle $H\subset TM$ is maximally non--integrable,
and hence defines a CR structure on $M$. From the construction it is
clear that a biholomorphism $f:\tilde M\to\tilde M$ which preserves
the hypersurface $M$ restricts to a CR automorphism of $M$.

The simplest example of this situation is provided by the unit sphere
$S^3\subset\Bbb C^2$. For $x\in S^3$ we get $T_xS^3=\{y\in\Bbb
C^2:\text{Re}(\langle x,y\rangle)=0\}$. The maximal complex subspace
of this is $H_x=\{y\in\Bbb C^2:\langle x,y\rangle=0\}$. One easily
verifies directly that this defines a contact structure on $S^3$. Hence
we have obtained a CR structure on $S^3$, called the \textit{standard
  structure}. CR structures which a locally isomorphic to the standard
structure on $S^3$ are called \textit{spherical}.

Any element $A\in U(2)$ defines a biholomorphism of $\Bbb C^2$ which
preserves the unit sphere $S^3$, and hence restricts to a CR
automorphism of the standard CR structure on $S^3$. Of course, this
action is transitive, so we see that $S^3$ with its standard structure
is a homogeneous CR manifold. The group of CR automorphisms of this
structure however is larger than $U(2)$. Identifying $S^3$ with the
space those complex lines in $\Bbb C^3$ which are isotropic for a
Hermitian inner product of signature $(2,1)$ leads to an faithful
action of $G:=PSU(2,1)$ on $S^3$ by CR automorphisms. Correspondingly,
one obtains a diffeomorphism $S^3\cong G/P$, where $P\subset G$ is the
stabilizer of an isotropic line in $\Bbb C^3$.

\subsection{Left invariant deformations of the standard structure}\label{2.2}
Restricting the action of $U(2)$ on $S^3$ further to $K:=SU(2)$ we
obtain a diffeomorphism $K\to S^3$, which we can use to carry over the
standard CR structure to $K$. In this picture, multiplication from the
left by any element of $K$ is a CR automorphism, so we have
constructed a left invariant CR structure on $K$. 

It is well known that left invariant structures on a Lie group can be
described in terms of the Lie algebra. Denoting by $e\in K$ the unit
element and by $\frak k=T_eK$ the Lie algebra of $K$, we get the fiber
$H_e\subset\frak k$ of the subbundle.  This must be a complex subspace
of complex dimension $1$ in the real vector space $\frak k$. By left
invariance, the fiber $H_g$ in each point $g\in K$ is spanned by the
values $L_X(g)$ of the left invariant vector fields generated by
elements $X\in H_e$, and the complex structure on $H_g$ comes from the
linear isomorphism $X\mapsto L_X(g)$. Explicitly, $\frak k$ consists
of all skew Hermitian $2\x 2$--matrices, i.e.
$$
\frak k=\left\{\begin{pmatrix} it & -\overline z\\ z &
    -it\end{pmatrix}:t\in\Bbb R,z\in\Bbb C\right\},
$$
and we will denote elements of $\frak k$ as pairs $(it,z)$. Using the
action on the first vector in the standard basis of $\Bbb C^2$ to
identify $K$ with $S^3$, we see that $H_e=\{(0,z):z\in\Bbb
C\}\subset\frak k$. The fact that this defines a left invariant
contact structure on $K$ is then immediate from the fact that
$[L_X,L_Y]=L_{[X,Y]}$ for all $X,Y\in\frak k$ and from
$[(0,1),(0,i)]=(-2i,0)$. Indeed, the linear functional $\al:\frak
k\to\Bbb R$ defined by $\al(it,z)=t$ defines a left invariant contact
form for the contact structure $H$. 

Now the crucial idea is that we can leave this left invariant contact
structure unchanged but deform the complex structure in the space
$H_e$ to obtain a family $(H,J_\la)$ of left invariant CR structures
on $K$ parametrized by a positive real number $\la$. Namely, for
$\la>0$ we define $J_\la(e)(0,u+iv):=(0,i(\la
u+i\tfrac{1}{\la}v))=(0,-\tfrac{1}{\la}v+i\la u)$. This extends to a
left invariant complex structure on the contact subbundle $H\subset
TK$, which in addition induces the standard orientation. The obvious
question is whether this is a true deformation of the standard CR
structure, or whether one just obtains (locally) isomorphic
structures.

Notice that, viewed as CR structures on $S^3$, the structures
$(H,J_\la)$ for $\la\neq 1$ are not invariant under the group $U(2)$.
Indeed the element $\left(\begin{smallmatrix} 1 & 0 \\ 0 &
    i\end{smallmatrix}\right)\in U(2)$ fixes the first vector in the
standard basis. The tangent map of its action is given by
$(it,z)\mapsto (it,iz)$, which is complex linear for $J_\la(e)$ if and
only if $\la=1$. Invariance under $U(2)$ would actually imply that the
structure is spherical, since by a classical result of Cartan, the
automorphism group of a non--spherical CR structure has dimension at
most three. A simple proof of this result can be found in
\cite{Srni04}.

\section{The canonical Cartan connections}\label{3}

\subsection{Three dimensional CR structures and Cartan
  geometries}\label{3.1}  
Three dimensional CR structures can be equivalently described as
Cartan geometries, which in particular implies that the curvature
gives a complete obstruction to being spherical. We first have to
describe the group $G=PSU(2,1)$ and its Lie algebra $\frak
g=\frak{su}(2,1)$ in a bit more detail. Consider the Hermitian form on
$\Bbb C^3$ defined by
$$
((z_0,z_1,z_2),(w_0,w_1,w_2))\mapsto
z_0\overline{w_2}+z_2\overline{w_0}+z_1\overline{w_1}. 
$$
Then the first and last vector in the standard basis are isotropic,
while the second one has positive length, so this form has signature
$(2,1)$. A direct computation shows that for this form we get
$$
\frak g=\left\{
  \begin{pmatrix}
  \al+i\be & w & i\ps \\ x & -2i\be & -\overline{w}\\ i\ph &
  -\overline{x} & -\al+i\be
  \end{pmatrix}: \al,\be,\ph,\ps\in\Bbb R, x,w\in\Bbb C
\right\}
$$
We obtain a grading $\frak g=\frak g_{-2}\oplus\frak g_{-1}\oplus\frak
g_0\oplus\frak g_1\oplus \frak g_2$ of $\frak g$ by
$$
\begin{pmatrix}
  \fg_0 & \fg_1 &\fg_2\\ \fg_{-1} & \fg_0 &\fg_1\\ \fg_{-2} &\fg_{-1}
  &\fg_0 
\end{pmatrix}.
$$
The associated filtration is defined by $\frak g^i=\frak
g_i\oplus\dots\oplus \frak g_2$, so we have 
$$
\frak g=\frak g^{-2}\supset\fg^{-1}\supset\dots\supset\frak g^2,
$$ 
and $[\fg^i,\fg^j]\subset\fg^{i+j}$. The parabolic subgroup $P\subset
G$ is the stabilizer of an isotropic line, for which we take the line
generated by the first basis vector. The Lie algebra of $P$ then
evidently is given by $\frak p=\fg^0=\fg_0\oplus\fg_1\oplus\fg_2$. In
particular the filtration $\{\fg^i\}$ is invariant under the adjoint
actions of $\frak p$ and $P$. 

\begin{definition*} (1) A \textit{Cartan geometry} of type $(G,P)$ on a
smooth manifold $M$ is a principal $P$--bundle $p:\Cal G\to M$
together with a one form $\om\in\Om^1(\Cal G,\frak g)$ such that
\begin{itemize}
\item $(r^g)^*\om=\Ad(g)^{-1}\o\om$ for all $g\in P$, where $r^g$
  denotes the principal right action of $g$.
\item $\om(\ze_A)=A$ for all $A\in\frak p$, where $\ze_A$ denotes the
  fundamental vector field with generator $A$.
\item $\om(u):T_u\Cal G\to\frak g$ is a linear isomorphism for all
$u\in\Cal G$.
\end{itemize}

\noindent
(2) A \textit{morphism} between two Cartan geometries $(\Cal G\to
M,\om)$ and $(\tcg\to\tilde M,\tilde\om)$ is a principal bundle
homomorphism $\Ph:\Cal G\to\tcg$ such that $\Ph^*\tilde\om=\om$. Note
that since both $\om$ and $\tilde\om$ are bijective on each tangent
space, this implies that $\Ph$ is a local diffeomorphism.

\noindent
(3) The \textit{homogeneous model} of the geometry is the principal
bundle $G\to G/P$ together with the left Maurer--Cartan form
$\om^{MC}$.
\end{definition*}

Given a Cartan geometry $(p:\Cal G\to M,\om)$ of type $(G,P)$ on $M$,
we can form the associated bundle $\Cal G\x_P(\frak g/\frak p)$. The
map $\Cal G\x (\frak g/\frak p)\to TM$ given by $(u,X)\mapsto
T_up\cdot(\om(u)^{-1}(X))$ descends to an isomorphism $\Cal
G\x_P(\frak g/\frak p)\cong TM$. Now $\frak g/\frak p$ contains the
$P$--invariant subspace $\frak g^{-1}/\frak p$, so this gives rise to
a subbundle $H\subset TM$. Moreover, $\frak g^{-1}/\frak p\cong\Bbb C$
and since $P$ consists of complex matrices, this complex structure is
invariant under the adjoint action of $P$. Therefore, it makes the
associated bundle $H=\Cal G\x_P(\frak g^{-1}/\frak p)$ into a complex
line bundle. If $H$ is a contact structure, then we obtain a three
dimensional CR structure on $M$.

\subsection{Regularity and normality}\label{3.2}
To characterize when $H$ is a contact structure, we need the curvature
$\ka\in\Om^2(\Cal G,\frak g)$ of $\om$. This is defined by
$\ka(\xi,\eta)=d\om(\xi,\eta)+[\om(\xi),\om(\eta)]$. In the case of
the homogeneous model, $\ka$ vanishes identically by the
Maurer--Cartan equation. Conversely, it can be shown that any Cartan
geometry with vanishing curvature is locally isomorphic to the
homogeneous model. Now we call a Cartan geometry of type $(G,P)$
\textit{regular} if and only if $\ka(\xi,\eta)$ has values in $\frak
g^{-1}\subset\frak g$ provided that $\om(\xi)$ and $\om(\eta)$ have
values in $\frak g^{-1}$.

Suppose that this condition is satisfied and that $\xi$ and $\eta$ are
local lifts of vector fields on $M$. Then the fact that $\om(\xi)$ and
$\om(\eta)$ have values in $\frak g^{-1}$ exactly means that these
vector fields are sections of $H\subset TM$. By definition of the
curvature and the assumptions we see that
$-\om([\xi,\eta])+[\om(\xi),\om(\eta)]$ has values in $\frak
g^{-1}$. Since $[\xi,\eta]$ lifts the bracket of the original fields,
this bracket cannot have values in $H$ unless $[\om(\xi),\om(\eta)]$
has values in $\frak g^{-1}$. One immediately verifies that the
bracket in $\frak g$ induces a non--degenerate map $\fg^{-1}/\frak p\x
\fg^{-1}/\frak p\to \fg/\fg^{-1}$. Hence we see that regularity of the
Cartan geometry ensures the we obtain an underlying CR structure. 

It is a general theorem, that any three dimensional CR structure
arises as the underlying structure of a Cartan geometry of type
$(G,P)$. However, there are many non--isomorphic Cartan geometries
having the same underlying CR structure. To get rid of this freedom,
one has to put an additional normalization condition on (the curvature
of) the Cartan connection $\om$. Under this additional condition, the
Cartan geometry is then uniquely determined up to isomorphism. See
\cite{Srni05} for a discussion of all these issues and
\cite{Cap-Schichl} for proofs, both in the realm of general parabolic
geometries.

We will not need the detailed form of the normalization condition, but
only some of its consequences. These follow from the fact that one may
relate the values of the curvature of a regular normal Cartan geometry
to certain explicitly computable Lie algebra cohomology groups. In the
case of three dimensional CR structures, these conditions imply that
$\ka(\xi,\eta)$ has values in $\frak g^1\subset\frak g$ for all $\xi$
and $\eta$. Moreover, if both $\om(\xi)$ and $\om(\eta)$ have values
in $\frak g^{-1}$, then $\ka(\xi,\eta)$ has to vanish identically.
Moreover, projecting the values of $\ka$ to $\frak g^1/\frak
g^2\cong\frak g_1$, one obtains the \textit{harmonic curvature}, which
still is a complete obstruction to the CR structure being spherical.

\subsection{The case of left invariant structures}\label{3.3}
Let us now consider one of the left invariant CR structures
$(H,J_\la)$ on $K=SU(2)$. As an ansatz, we use the trivial principal
$P$--bundle $\Cal G:=K\x P$. For $X\in\frak k$, define $\hat
L_X:=(L_X,0)\in\frak X(\Cal G)$.  The second part of the ansatz is
that $\om(\hat L_X)$ is constant along $K\x\{e\}$. The motivation for
this ansatz is as follows. For any $k\in K$, the left translation by
$k$ defines a CR automorphism of $K$, which leaves each $L_X$
invariant.  These automorphisms lift to the canonical principal bundle
in a way compatible with the canonical Cartan connection. Fixing an
identification of the fiber of the Cartan bundle over $e\in K$ with
$P$, we can use these lifts to trivialize the Cartan bundle and in
such a way that $\om(\hat L_X)$ is constant along $K\x\{e\}$.

Consider a linear map $\ph:\frak k\to\frak g$ such that the
composition with the projection $\frak g\to\frak g/\frak p$ with $\ph$
is a linear isomorphism. Any tangent vector in $(k,g)\in K\x P$ can be
uniquely written as $(L_X(k),L_A(g))$ for some $X\in\frak k$ and
$A\in\frak p$. Hence we can define 
$\om\in\Om^1(K\x P,\frak g)$ by
$$
\om(L_X(k),L_A(g)):=\Ad(g^{-1})(\ph(X))+A.
$$
By the assumption on $\ph$ this defines a linear isomorphism on each
tangent space, and using that the principal right action is just
multiplication from the right in the second factor and that
$\ze_A=(0,L_A)$ for each $A\in\frak p$, one immediately verifies that
this defines a Cartan connection. 

We can also immediately compute the curvature $\ka$ of this
connection. Since $\ka$ is horizontal and $P$--equivariant, it
suffices to compute $\ka(\hat L_X,\hat L_Y)(k,e)$ for all $X,Y\in\frak
k$ and $k\in K$. Now by definition,
$$
\ka(\hat L_X,\hat L_Y)=d\om(\hat L_X,\hat L_Y)+[\om(\hat
L_X),\om(\hat L_Y)]=-\om([\hat L_X,\hat L_Y])+ [\om(\hat L_X),\om(\hat
L_Y)].
$$
Using $[L_X,L_Y]=L_{[X,Y]}$ we see that, along $K\x\{e\}$, the
function $\ka(\hat L_X,\hat L_Y)$ is constant and equal to
$$
[\ph(X),\ph(Y)]-\ph([X,Y]).
$$
Hence the curvature exactly expresses the obstruction against $\ph$
being a homomorphism of Lie algebras. 

It remains to express the fact that the Cartan connection $\om$
induces the ``right'' underlying CR structure in terms of the linear
map $\ph$. Returning to the notation of \ref{2.2}, we denote elements
$X\in\frak k$ as pairs $(it,z)$ for $t\in\Bbb R$ and $z\in\Bbb C$.
Then $L_{(it,z)}$ lies in the contact subbundle $H$ if and only if
$t=0$, so to get the right contact subbundle, $\ph(it,z)$ must lie in
the subspace $\frak g^{-1}\subset\frak g$ if and only if $t=0$.  Given
this, we get an induced linear map $\frak k\supset H_e\to\frak
g^{-1}$. Composing with the natural projection, we get a linear
isomorphism $\frak k\to\frak g^{-1}/\frak p\cong\Bbb C$. The condition
that we get the induced complex structure $J_\la$ exactly means that
via this isomorphism the (fixed) standard complex structure on $\frak
g^{-1}/\frak p$ induces the complex structure $J_\la(e)$ on $H_e$.

Having all this at hand, we can prove the main technical result of
this article:
\begin{thm*}
  For fixed $\la>0$, the linear map $\ph_\la:\frak k\to\frak g$
  defined by 
$$
\ph_\la(it,u+iv):=\begin{pmatrix} \frac{1+\la^2}{4\la} it &
  -\frac{5-3\la^2}{4\sqrt{\la}}u-\frac{3-5\la^2}{4\la\sqrt{\la}}iv &
  \frac{-15+34\la^2-15\la^4}{16\la^2} it\\  
  \sqrt{\la} u+\frac{1}{\sqrt{\la}}iv & -\frac{1+\la^2}{2\la} it &
  \frac{5-3\la^2}{4\sqrt{\la}}u-\frac{3-5\la^2}{4\la\sqrt{\la}}iv \\
  it & -\sqrt{\la} u+\frac{1}{\sqrt{\la}}iv &  \frac{1+\la^2}{4\la} it\end{pmatrix}
$$
induces a linear isomorphism $\frak k\to\frak g/\frak p$. It has
the property that $\ph_\la(H_e)\subset\fg^{-1}$, and via the
induced isomorphism $H_e\to\fg^{-1}/\frak p$ the induced complex
structure on $H_e$ is $J_\la(e)$. Finally, the map $\ka_\la:\frak
k\x\frak k\to\frak g$ defined by
$$
\ka_\la(X,Y):=[\ph_\la(X),\ph_\la(Y)]-\ph_\la([X,Y])
$$ 
has values in $\frak g^1$ and vanishes on $H_e\x H_e$.

Explicitly, $\ka_\la$ satisfies 
$$
\ka_\la ((it,0),(0,u+iv))=\begin{pmatrix} 0 &
  -\frac{3t(\la^4-1)}{2\la^2\sqrt{\la}}(v-i\la u) & 0 \\
  0 & 0 & \frac{3t(\la^4-1)}{2\la^2\sqrt{\la}}(v+i\la u)\\ 0 & 0 &
  0\end{pmatrix}
$$
and this completely determines $\ka_\la$.
\end{thm*}
\begin{proof}
  From the definition of $\ph_\la$ it is evident that it induces a
  linear isomorphism $\frak k\to\frak g/\frak p$ and that it maps
  elements of $H_e$, which are characterized by $it=0$ to $\frak
  g^{-1}$. Then the isomorphism $H_e\to \frak g^{-1}/\frak p$ is given
  by $u+iv\mapsto \sqrt{\la} u+i\frac{1}{\sqrt{\la}}v$, so the complex
  structure on $\fg^{-1}/\frak p$ evidently induces $J_\la(e)$ on
  $H_e$. It is then straightforward but tedious to check that
  $\ka_\la$ has values in $\frak g^1$ and vanishes on $H_e\x H_e$ as
  well as the explicit formula. That this expression determines
  $\ka_\la$ follows since by skew symmetry and vanishing of $\ka_\la$
  on $H_e\x H_e$ we obtain
$$
\ka_\la((it,z),(it',z'))=\ka_\la((it,0),(0,z'))-\ka_\la((it',0),(0,z)). 
$$
\end{proof}

\subsection{Digression: How to get the formula for
  $\ph_\la$}\label{3.4} 

The result of Theorem \ref{3.3} is all that is needed in the sequel.
Since the proof does not explain how the formula for $\ph_\la$ was
obtained (although really doing the computation gives some hints), we
will briefly discuss this. As a spin off, this will show that
$\ph_\la$ is essentially uniquely determined by the four properties
listed in Theorem \ref{3.3}. The main point is that there is some
evident non--uniqueness around, and dealing with this is the key step
to determine $\ph_\la$. Recall that for any element $g\in P$, the
adjoint action $\Ad(g):\fg\to\fg$ preserves the filtration. In
particular, it preserves $\frak g^1$ and $\frak g^{-1}$ as well as
$\frak p$ and therefore induces a linear isomorphisms on $\frak
g/\frak p$ and $\frak g^{-1}/\frak p$. One immediately checks that the
second of these isomorphisms is complex linear. From these
observations it is evident, that if $\ph:\frak k\to\frak g$ is a
linear map which satisfies the four properties of Theorem \ref{3.3}
and $g\in P$ is arbitrary, then also $\Ad(g)\o\ph$ has these
properties. 

To deal with this freedom, we need a bit more information on the group
$P$. Note first that there is a subgroup $G_0\subset P$ consisting of
all $g\in P$ for which $\Ad(g):\fg\to\fg$ even preserves the
grading. It is a general result (see \cite[Proposition
2.10]{Cap-Schichl}) that $G_0$ has Lie algebra $\frak g_0$ and any
$g\in G$ can be uniquely written in the form $g=g_0\exp(Z_1)\exp(Z_2)$
for $g_0\in G_0$ and $Z_i\in\frak g_i$. For our choice of $G$ and $P$,
one immediately verifies that the (complex) linear automorphism on
$\frak g^{-1}/\frak p$ induced by $\Ad(g)$ depends only on $g_0$ and
one obtains an isomorphism $G_0\to\Bbb C\setminus\{0\}$ in this way. 

But now any linear isomorphism $H_e\to\frak g^{-1}/\frak p$,
which induces $J_\la(e)$ on $H_e$ can be written (identifying $
\frak g^{-1}/\frak p$ with the matrix component in the first column of
the second row) as the composition of a complex linear automorphism of
$\frak g^{-1}/\frak p$ with $u+iv\mapsto \sqrt{\la}
u+\frac{1}{\sqrt{\la}}iv$.

Hence if we want $\ph$ to induce a linear isomorphism $\frak k\to\frak
g/\frak p$, map $H_e\to\frak g^{-1}$, and induce $J_\la(e)$, then
we may assume the the lower two rows of the first column of
$\ph(it,u+iv)$ have the form $\begin{pmatrix} \sqrt{\la}
  u+\frac{1}{\sqrt{\la}}iv+tz_0\\ is t\end{pmatrix}$ for some
$z_0\in\Bbb C$ and some $s\in\Bbb R\setminus\{0\}$. (Of course, this
also determines the second component in the last row.) Making this
ansatz also reduces the freedom to composition with
$\Ad(\exp(Z_1)\exp(Z_2))$. Having made this ansatz, one can already
compute the $\frak g_{-2}$--component of the restriction of $\ka$ to
$H_e\x H_e$ and vanishing of this forces $s=1$.

Next we observe taking the bracket with a nonzero element of $\frak
g_{-2}$ induces a linear isomorphism $\frak g_1\to\fg_{-1}$. Using
this, we see that composing with $\Ad(\exp(Z_1))$ for an appropriate
choice of $Z_1$ we can require $z_0=0$ in the above ansatz, and this
reduces the freedom to composition with $\Ad(\exp(Z_2))$. To get rid
of this freedom, we observe that bracketing with a nonzero element of
$\frak g_{-2}$ induces a linear isomorphism from $\frak g_2$ to the
(one--dimensional) space of real diagonal matrices contained in $\frak
g$. Hence we can eliminate all the freedom of composition with $\Ad(g)$
by the ansatz that the first column of $\ph(it,u+iv)$ has the form
$\begin{pmatrix} uz_0+vz_1+ist\\ \sqrt{\la} u+\frac{1}{\sqrt{\la}}iv\\
  it\end{pmatrix}$ for elements $z_0,z_1\in\Bbb C$ and $s\in\Bbb R$.
Having made this ansatz, one can compute the complete $\frak g_{-2}$
component of $\ka$ and the $\frak g_{-1}$ component of the restriction
to $H_e\x H_e$, and vanishing of these forces $z_0=z_1=0$.

Now one can, step by step, take ansatzes for the remaining components
of $\ph(it,u+iv)$ and determine components of $\ka$. In the end, one
finds out that the conditions on $\ka$ in Theorem \ref{3.3} are
sufficient to uniquely pin down the formula for $\ph_\la$.

\subsection{The canonical Cartan connections}\label{3.5}
It is now easy to show that the map $\ph_\la$ from Theorem \ref{3.3}
leads to the canonical Cartan connection for $(K,H,J_\la)$.

\begin{kor*}
  (1) For some $\la>0$, consider the left invariant CR structure
  $(H,J_\la)$ on $K=SU(2)$ from \ref{2.2}. Then the regular normal
  parabolic geometry associated to this structure is $(K\x P\to
  K,\om_\la)$, where
$$
\om_\la(L_X(k),L_A(g))=\Ad(g^{-1})(\ph_\la(X))+A
$$
with $\ph_\la:\frak k\to\frak g$ the map from Theorem \ref{3.3}.

\noindent
(2) The CR structure $(H,J_\la)$ is spherical if and only if
    $\la=1$, i.e.~if and only if it equals the standard structure.  
\end{kor*}
\begin{proof}
  (1) From \ref{3.3} we know that the formula for $\om_\la$ defines a
  Cartan connection on the trivial bundle $K\x P$. The conditions on
  $\ph_\la$ in Theorem \ref{3.3} which do not involve $\ka_\la$
  exactly say the this Cartan connection induces the CR structure
  $(H,J_\la)$ on $K$. Hence to prove (1), it remains to show that
  $\om_\la$ is normal. The formula for $\ka_\la$ in Theorem \ref{3.3}
  gives us the restriction of the curvature of $\om_\la$ to
  $K\x\{e\}$. By equivariancy of the normalization condition it
  suffices to show normality of this restriction in order to prove
  that $\om_\la$ is normal. Since $\ka_\la$ has values in $\fg^1$ and
  the restriction to $H_e\x H_e$ vanishes, it is homogeneous of degree
  $\geq 4$, and the component of degree $4$ maps $(\frak k/H_e)\x H_e$
  to $\frak g_1$. Identifying $\frak g_1$ with the component in the
  second column of the first row of a matrix, this component is
  complex linear in the second variable (with respect to $J_\la$). It
  is well known (and easy to see) that the one dimensional space of
  such maps exactly constitutes the harmonic part in degree $4$, so in
  particular such maps lie in the kernel of the Kostant
  codifferential. Since maps of homogeneity $\geq 5$ automatically
  have that property, normality follows.

\noindent
(2) This is now evident since $\ka_\la$ vanishes if and only if
    $\la^4=1$. 
\end{proof}

Notice that we can use the same construction replacing $G$ by the
three--fold covering $SU(2,1)$ and $P$ by the stabilizer of a line in
that group. Such an extension is necessary for example if one wants to
have a standard tractor bundle, compare with \cite{Gover-Graham}. In
the picture of CR geometry, such an extension is associated to the
choice of a third root of a certain complex line bundle. In our case,
this bundle is trivial, so this poses no problem.

\subsection{Tractors and tractor calculus}\label{3.6}
As an indication how the description of the canonical Cartan
connection in Corollary \ref{3.5} can be used further, we discuss
tractor bundles and compute normal tractor connections. We will work
here in the setting that $G=SU(2,1)$ and $P\subset G$ is the
stabilizer of a line. Recall that for a representation $V$ of the
group $G$, one obtains a tractor bundle by restricting the
representation to $P\subset G$ and forming the associated bundle to
the canonical Cartan bundle. While sections of these bundles are
unusual geometric objects, they have the advantage that they carry
canonical linear connections induced by the canonical Cartan
connection. 
\begin{prop*}
  For some $\la>0$ consider the left invariant CR structure
  $(H,J_\la)$ on $K=SU(2)$ from \ref{2.2}, and let $V$ be a
  representation of $G=SU(2,1)$. Then the associated tractor bundle
  $\Cal T\to K$ is canonically trivial, so $\Ga(\Cal T)\cong
  C^\infty(K,V)$. In this identification, the tractor connection
  $\nabla^{\Cal T}$ is determined by
$$
\nabla^{\Cal T}_{L_X}f=L_X\cdot f+\rho(\ph_\la(X))\o f, 
$$ 
for $f:K\to V$, where $L_X\in\frak X(K)$ denotes the left invariant
vector field generated by $X\in\frak k$ and $\rho:\frak g\to L(V,V)$
is the derivative of the representation of $G$ on $V$.
\end{prop*}
\begin{proof}
Since the canonical Cartan bundle is trivial, so is the associated
bundle $\Cal T$. Explicitly, the identification $\Ga(\Cal T)\to
C^\infty(K,V)$ is given by restricting the $P$--equivariant function
$\Cal G=K\x P\to V$ corresponding to a section to the subset
$K\x\{e\}$. In terms of equivariant functions, the tractor connection
can be easily described explicitly, see \cite[section 3]{tractors}:
For the equivariant map $h:\Cal G\to V$ corresponding to $s\in\Ga(\Cal
T)$, a vector field $\xi$ on $K$ and a lift $\tilde\xi\in\frak X(\Cal
G)$ of $\xi$, the covariant derivative $\nabla^{\Cal T}_\xi s$ is
represented by the function $\tilde\xi\cdot h+\rho(\om(\xi))\o h$.

Putting $\xi=L_X$, we can use $(L_X,0)$ for $\tilde\xi$. This has the
particular advantage that its flow leaves the subset $K\x\{0\}\subset
K\x P$ invariant. Therefore, putting $f:=h|_{K\x\{e\}}$ we see that
$((L_X,0)\cdot h)|_{K\x\{e\}}=L_X\cdot f$.For the second term,
restriction to $K\x\{e\}$ makes no problems anyhow, so the formula for
$\nabla^{\Cal T}$ follows.
\end{proof}

As a concrete example, let us describe how the three dimensional
family of infinitesimal automorphisms corresponding to the left
translations by elements of $K$ are represented within adjoint
tractors. This means that we consider the representation $V=\frak g$,
and the resulting tractor bundle is the adjoint tractor bundle $\Cal
A$. The canonical Cartan connection induces an isomorphism between
$\Ga(\Cal A)$ and the space of right invariant vector fields on $\Cal
G$, see \cite{deformations}. Infinitesimal automorphisms of a Cartan
geometry are described by such vector fields, and they are
characterized by a simple differential equation, see \cite[Proposition
3.2]{deformations}. We will verify this equation for the three
dimensional family of infinitesimal automorphisms corresponding to
left translations on $K$.

The construction of the canonical Cartan connection on $\Cal G=K\x P$
for the left invariant CR structure $(H,J_\la)$ on $K$ shows that for
each $k'\in K$ the map $(k,g)\mapsto (k'k,g)$ is the lift of the left
translation by $k'$ to an automorphism of the parabolic geometry
$(\Cal G,\om_{\la})$. The infinitesimal generators of this three
parameter group of automorphisms are of course the vector fields
$(R_X,0)$ for $X\in\frak k$, where $R_X$ denotes the right invariant
vector field. Let $s_X\in\Ga(\Cal A)$ be the corresponding section,
i.e.~the smooth equivariant function corresponding to $s_X$ is
$\om_\la((R_X,0))$. Since $R_X(k)=L_{\Ad(k^{-1})X}(k)$ we see that the
smooth function $f_X:K\to\frak g$ corresponding to $s_X$ is given by
$f_X(k)=\ph_\la(\Ad(k^{-1})X)$. From the proposition above, we
conclude that $\nabla^{\Cal A}_{L_Y}s_X$ corresponds to the function
$$
L_Y\cdot\ph_\la(\Ad(k^{-1})X)+[\ph_\la(Y),\ph_\la(\Ad(k^{-1})X)].
$$
Now the first term can be computed as
$$
\tfrac{d}{dt}|_{t=0}\ph_\la(\Ad(\exp(-tY))\Ad(k^{-1})X)=-\ph_\la([Y,\Ad(k^{-1})X]).
$$
Hence we see that $\nabla^{\Cal A}_{L_Y}s_X$ corresponds to the
function $\ka_\la(Y,\Ad(k^{-1})(X))$ which represents the curvature of
$\om_\la$ evaluated on the vector fields $L_Y$ and $R_X$. This is
exactly the infinitesimal automorphism equation from \cite[Proposition
  3.2]{deformations}.

\end{document}